\documentclass[10pt]{article}
\usepackage[all]{xy}
\usepackage{amsfonts,amsmath,oldgerm,amssymb,amscd}
\newcommand{\ra}{\rightarrow}		
\newcommand{\by}[1]{\stackrel{#1}{\ra}}

\newcommand{\surj}{\ra\!\!\!\ra}	
\newcommand{\ol}{\overline}		\newcommand{\wt}{\widetilde}
\newcommand{\iso}{\by \sim}

\newtheorem{theorem}{Theorem}[section]
\newtheorem{proposition}[theorem]{Proposition}
\newtheorem{lemma}[theorem]{Lemma}

\newtheorem{corollary}[theorem]{Corollary}

\newcommand{\ga}{\alpha}	
		
\newcommand{\gf}{\varphi}	
	\newcommand{\gl}{\lambda}

	\newcommand{\BN}{\mbox{$\mathbb N$}}
	
	\newcommand{\BR}{\mbox{$\mathbb R$}}

	\newcommand{\BZ}{\mbox{$\mathbb Z$}}

\newcommand{\CC}{\mbox{$\mathcal C$}}

\newcommand{\mm}{\mbox{$\mathfrak m$}}	
	\newcommand{\p}{\mbox{$\mathfrak p$}}

\newcommand{\ot}{\mbox{\,$\otimes$\,}}	\newcommand{\op}{\,\mbox{$\oplus$\,}}
\newcommand{\Spec}{\mbox{\rm Spec\,}}	\newcommand{\hh}{\mbox{\rm ht\,}}
	
\newcommand{\Aut}{\mbox{\rm Aut}}	
\newcommand{\Hom}{\mbox{\rm Hom}}	\newcommand{\Sing}{\mbox{\rm Sing\,}}

\newcommand{\Um}{\mbox{\rm Um}}		\newcommand{\SL}{\mbox{\rm SL}}
\newcommand{\GL}{\mbox{\rm GL}}		
\newcommand{\tr}{\mbox{\small \BR}}

\def\remark{\refstepcounter{theorem}\paragraph{Remark \thetheorem}}

\begin{document}

\begin{center}

{\Large A note on projective modules over
real affine algebras} \\
\vspace{.2in}

{\large Manoj Kumar Keshari} 
\footnote{Partially supported by Kanwal Rekhi career development 
scholarship of the TIFR endowment fund.}   

\vspace{.1in}
{\small \it School of Mathematics, Tata Institute of Fundamental Research, \\
Homi Bhabha Road, Mumbai-400005, India. \\
e-mail :  manoj@math.tifr.res.in}
\end{center}          

\section{Introduction}
Let $A$ be an affine domain over $\BR$ of
dimension $d$. Let $f\in A$ be an element not belonging to any real
maximal ideal of $A$ and let
$P$ be a projective $A$-module of rank $\geq d-1$. Let
$(a,p)\in A_f\op P_f$ be a unimodular element and $Q=A_f\op
P_f/(a,p)A_f$. If $P$ is free, then a result of Ojanguren and Parimala
($\cite{O-P}$, Theorem) shows that $Q$ is extended from $A$. A
consequence of this result is that, if $d=3$, then all projective modules over
$A_f$ are free, where $A=\BR[X_1,X_2,X_3]$ (see $\cite{O-P}$ for
motivation). In this paper, we prove the following result (\ref{cor})
which is a generalization of the above result of Ojanguren and Parimala.
\\

\noindent {\bf Theorem :}
{\it Let $A$ be an affine algebra over
$\BR$ of dimension $d$. Let $f\in A$ be an element not belonging to
any real maximal ideal of $A$. Let $P$ be a projective $A$-module
of rank $\geq d-1$. Let $(a,p)\in A_{f}\op P_f$ be a
unimodular element. Then, the projective $A_{f}$-module $Q=A_{f}\op
P_f/(a,p)A_{f}$ is extended from $A$. }


\section{Preliminaries}

In this paper, all the rings are assumed to be commutative Noetherian
and all the modules are finitely generated. 

Let $B$ be a ring and let $P$ be a projective $B$-module. Recall that $p\in
P$ is called a unimodular element if there exists an $\psi \in
P^*=\Hom_B(P,B)$ such that $\psi(p)=1$. We denote by $\Um(P)$, the set
of all unimodular elements of $P$. We write $O(p)$ for the ideal of
$B$ generated by $\psi(p)$, for all $\psi \in P^*$. Note that, if
$p\in P$ is a unimodular element, then $O(p)=B$.
 
Let $E_n(B)$ denote the subgroup of $\SL_n(B)$ generated by
all the elementary matrices $E_{ij}(z)$, where
$E_{ij}(z)\in \SL_n(B)$ is such that its diagonal elements are
$1$, $i\neq j$, $(i,j)$th entry is $z$ and the rest of the entries are $0$,
where $z\in B$. 
\medskip

We begin by  stating a classical result of Serre \cite{S1}.

\begin{theorem}\label{serre}
Let $A$ be a ring of dimension $d$. Then, any projective
$A$-module $P$ of rank $> d$ has a unimodular element. In particular,
if $\dim A=1$, then any projective $A$-module of trivial determinant is
free.
\end{theorem}

Let $B$ be a ring and let $P$ be a projective $B$-module. 
Given an element
$\gf\in P^\ast$ and an element $p\in P$, we define an endomorphism
$\gf_p$ of $P$ as the composite $P\by \gf B\by p P$. 

If $\gf(p)=0$, then ${\gf_p}^2=0$ and, hence, $1+\gf_p$ is a unipotent 
automorphism of $P$.

By a ``transvection'', we mean an automorphism of $P$ of the form
$1+\gf_p$, where $\gf(p)=0$ and either $\gf$ is unimodular in $P^\ast$
or $p$ is unimodular in $P$. We denote by $E(P)$, the subgroup of
$\Aut(P)$ generated by all transvections of $P$. Note that, $E(P)$ is a
normal subgroup of $\Aut(P)$. 

An existence of a transvection of $P$ pre-supposes that $P$ has a
unimodular element. Now, let $P = B\op Q$, $q\in Q, \alpha\in
Q^*$. Then $\Delta_q(b,q')=(b,q'+bq)$ and
$\Gamma_\alpha(b,q')=(b+\alpha(q'),q')$ are transvections of
$P$. Conversely, any transvection $\Theta$ of $P$ gives rise to a
decomposition  $P=B\op Q$ in such a way that $\Theta = \Delta_q$ or
$\Theta = \Gamma_\alpha$.
\medskip

Now, we state a classical result of Bass $\cite{B}$.

\begin{theorem}\label{bass}
Let $A$ be a ring of dimension $d$ and let $P$ be a
projective $A$-module of rank $> d$. Then $E(A\op P)$ acts
transitively on $\Um(A\op P)$. 
\end{theorem}

The following result is due to Bhatwadekar and Roy
($\cite{Bh-Roy}$, Proposition 4.1) and is 
about lifting an automorphism of a projective module. 

\begin{proposition}\label{one}
Let $A$ be a ring and let $J$ be an ideal of $A$. Let $P$ be a
projective $A$-module of rank $n$. 
Then, any transvection $\wt \Theta$ of $P/JP$ (i.e. $\wt \Theta \in E(P/JP)$)
can be lifted to a (unipotent) automorphism $\Theta$ of $P$. In
particular, if $P/JP$ is free (of rank $n$), then any element $\ol \Psi$ of
$E((A/J)^n)$ can be lifted to $\Psi \in \Aut(P)$. If in
addition, the natural map $\Um(P) \ra \Um(P/JP)$ is surjective, then
the natural map $E(P)\ra E(P/JP)$ is surjective.
\end{proposition}

Now, we recall some preliminary facts about symplectic modules. 
Let $A$ be a ring. A bilinear map $\langle ,\rangle : A^n\times A^n
\ra A$ is called {\it alternating} if $\langle v,v\rangle=0$,
$\forall\, v\in A^n$. Let us fix a basis of $A^n$, say
$e_1,\ldots,e_n$. Let $\langle e_i,e_j\rangle=a_{ij} \in A$. Then
$\ga=(a_{ij}) \in M_n(A)$ is such that $\ga+\ga^t=0$.  Thus, giving an
alternating bilinear form $\langle,\rangle$ on $A^n$ is equivalent to
giving a $n\times n$ matrix $\ga$ such that $\ga +
\ga^t=0$. Conversely, if $2\in A^*$ (the set of units of $A$),
 then an $n\times n$ matrix $\ga=(a_{ij})$
such that $\ga + \ga^t=0$ gives rise to a bilinear alternating map
$\langle,\rangle : A^n \times A^n \ra A$ given by
$\langle e_i,e_j\rangle =a_{ij}$.

An alternating form $\langle,\rangle$ on $A^n$ is called {\it
non-degenerate} if the corresponding $n\times n$ matrix $\ga$ is
invertible. A {\it Symplectic} $A$-module of rank $n$ is a pair
$(A^n,\langle,\rangle)$, where $\langle,\rangle : A^n \times A^n \ra
A$ is a non-degenerate alternating bilinear form. If
$(A^n,\langle,\rangle)$ is a symplectic $A$-module, then, it is easy 
to see that $n$ is even.

Two symplectic modules $(A^n,\langle,\rangle_1)$ and
$(A^n,\langle,\rangle_2)$ are said to be isomorphic if there exists an
isomorphism $\tau:A^n \iso A^n$ such that $\langle v_1,v_2\rangle_1 =
\langle\tau(v_1),\tau(v_2)\rangle_2$, $\forall \; v_1,v_2\in A^n$.

To make the notation simple, we will always denote a non-degenerate
alternating bilinear form  by $\langle,\rangle$. 

If $(A^n,\langle,\rangle)$ and $(A^m,\langle,\rangle)$ are two
symplectic modules, then non-degenerate alternating bilinear forms on
$A^n$ and $A^m$ will give rise (in a canonical manner) to a
non-degenerate alternating bilinear form on $A^n\op A^m =A^{n+m}$
 and we denote
the symplectic module thus obtained by $(A^n \perp
A^m,\langle,\rangle)$. There is a unique (upto scalar multiplication
by elements of $A^*$) non-degenerate alternating bilinear form
$\langle,\rangle$ on $A^2$, namely $\langle(a,b),(c,d)\rangle=ad-bc$.

An {\it isometry} of the symplectic module $(A^n,\langle,\rangle)$ is an
automorphism of $(A^n,\langle,\rangle)$. We denote by
$Sp_n(A,\langle,\rangle)$ the group of isometries of
$(A^n,\langle,\rangle)$. It is easy to see that
$Sp_n(A,\langle,\rangle)$ is a subgroup of $\SL_n(A)$ and it coincides
with $\SL_n(A)$ when $n=2$. Therefore, $\SL_2(A)$ can be identified
with a subgroup of $Sp(A^2\perp A^n,\langle,\rangle)$.

Let $(A^n,\langle,\rangle)$ be a symplectic $A$-module and let $u,v\in
A^n$ be such that $\langle u,v\rangle =0$. Let $a\in A$ and let
$\tau_{(a,u,v)} : A^n \ra A^n$ be a map defined by
$$\tau_{(a,u,v)} (w)=w + \langle w,v\rangle u+\langle w,u\rangle v+a
\langle w,u\rangle u, \;{\rm for}\; w\in A^n.$$ Then $\tau_{(a,u,v)} \in
Sp_n(A,\langle,\rangle)$. Moreover, it is easy to see that
$$\tau_{(a,u,v)}^{-1} = \tau_{(-a,-u,v)} = \tau_{(-a,u,-v)}\;\; {\mbox
{\rm and}}\;\; \ga \tau_{(a,u,v)}\ga^{-1} = \tau_{(a,\ga(u),\ga(v))}$$ for
an element $\ga \in Sp_n(A,\langle,\rangle)$.

An isometry $\tau_{(a,u,v)}$ is called a {\it symplectic transvection}
if either $u$ or $v$ is a unimodular element in $A^n$. We denote by
$ESp_n(A,\langle,\rangle)$ the subgroup of $Sp_n(A,\langle,\rangle)$
generated by symplectic transvections. It follows from the above
discussion that $ESp_n(A,\langle,\rangle)$ is a normal subgroup of
$Sp_n(A,\langle,\rangle)$.\\

The following result is due to Bhatwadekar ($\cite{Bh}$, Corollary 3.3) and
is about lifting an automorphism of a projective module. It is a
generalization of a result of Suslin ($\cite{Su2}$, Lemma 2.1).

\begin{proposition}\label{two}
Let $B$ be a two dimensional ring and let $I$ be an ideal of $B$ such
that $\dim (B/I)\leq 1$. Let $P$ be a projective $B$-module of
(constant) rank $2$ such that $P/IP$ is free. Then, any element of
$\SL_2(B/I) \cap ESp_4(B/I)$ can be lifted to an element of $\SL(P)$.
\end{proposition}

Let $A$ be a commutative ring and let $I$ be an ideal
 of $A$. For $n\geq 3$, let $E^1_n(A,I)$ denote the
subgroup of $E_n(A)$ generated by elementary matrices
$E_{1i}(a)$ and $E_{j1}(x)$, where
$2\leq i,j \leq n$, $a\in A$, $x\in I$. 

Let $\GL_n(A,I)$ denote the kernel of the canonical map $\GL_n(A) \ra
\GL_n(A/I)$. For $n\geq 3$, we denote $E_n^1(A,I)\cap \GL_n(A,I)$ by
$E_n(A,I)$. 

Let $P$ be a finitely generated
projective $A$-module of (constant) rank $d$. Let $t$ be a non zero
divisor of $A$ such that $P_t$ is free. Then, it is easy to see that
there exits a free submodule $F=A^d$ of $P$ and a positive integer $l$
such that, if $s=t^l$, then $sP\subset F$. Therefore, $sF^*\subset
P^*\subset F^*$.

\begin{lemma}\label{lift}
Let $A,P,F,s$ be as above. If $p\in F$, then $\Delta_p\in E(A\op F)
\cap E(A\op P)$ and if $\alpha \in F^*$, then $\Gamma_{s\alpha}\in
E(A\op F)\cap E(A\op P)$. Hence, if $d\geq 2$ and, if we identify
$E_{d+1}(A)$ with $E(A\op A^d)$, then $E_{d+1}^1(A,As)$ can be regarded
as a subgroup of $E(A\op P)$.
\end{lemma}

We denote by $\Um_n(A,I)$, the set of $I$-unimodular rows of length
$n$ over $A$ (i.e. unimodular rows of the type $[a_1,\ldots,a_n]$,
$1-a_1\in I$ and $a_i\in I, 2\leq i\leq n$).

For $n\geq 3$, $MSE_n(A,I)$ will denote the orbit set
$\Um_n(A,I)/E_n(A,I)$. We write $MSE_n(A)$ for $MSE_n(A,A)$.

Let $A$ be a commutative ring and let $I$ be an ideal of $A$. Let
$B=\BZ \op I$ (with the obvious ring structure on $B$). Then, for
$n\geq 3$, the canonical ring homomorphism $B\ra A$ gives rise to a
map $E_n(B,I)\ra E_n(A,I)$, a surjective map $\Um_n(B,I) \surj
\Um_n(A,I)$ and, hence, a surjective map $MSE_n(B,I) \surj MSE_n(A,I)$. 
\medskip

The following theorem is due to W. van der Kallen ($\cite{K}$, Theorem
3.21) and is very crucial for our result.

\begin{theorem}\label{kallen}
(Excision) Let $n\geq 3$. Let $A$ be a commutative ring and let $I$ be
an ideal of $A$. Then, the canonical maps $MSE_n(\BZ\op I,I) \ra
MSE_n(A,I)$ and $MSE_n(\BZ \op I,I) \ra MSE_n(\BZ\op I)$ are bijective.
\end{theorem}

The following result is due to Vaserstein ($\cite{V}$, Theorem).

\begin{theorem}\label{vaserstein}
Let $B$ be a commutative ring and let $[b_1,\ldots,b_n] \in
\Um_n(B)$, $n\geq 3$. Let $d$ be a positive integer. Then 
$$[b_1^d,b_2,\ldots,b_n] = [b_1,b_2^d,\ldots,b_n]\; (\rm{mod}\;
E_n(B)). $$
\end{theorem}

The following corollary is a consequence of Theorems $\ref{kallen}$ and
$\ref{vaserstein}$.  

\begin{corollary}
Let $A$ be a ring and  $I$  an ideal of $A$. Let
$[a_1,\ldots,a_n]$ be an element of $\Um_n(A,I), n\geq 3$. Let $d$ be
a positive integer. Then
$$[a_1^d,a_2,\ldots,a_n] = [a_1,a_2^d,\ldots,a_n]\; (\rm{mod}\;
E_n(A,I)).$$  
\end{corollary}

The following result of Suslin ($\cite{su1}$, Lemma 2) is also used in
the proof of our result $(\ref{main})$.

\begin{proposition}\label{su}
Let $A$ be a commutative ring and let $P$ be a finitely generated
projective $A$-module of rank $d$. Let $(c,p_1)\in A\op P$ be a
unimodular element. Suppose that $P/cP$ is a free $A/Ac$-module of
rank $d$ and that $\ol p_1 \in P/cP$ can be extended to a basis of
$P/cP$. Then, there exists an $A$-automorphism $\Phi$ of $A\op P$ such
that $\Phi (c^d,p_1)=(1,0)$.
\end{proposition}

The following result is due to Ojanguren and Parimala ($\cite{O-P}$,
Proposition 3). 

\begin{proposition}\label{su1}
Let $\CC = \Spec C$ be a smooth affine curve over a field $k$ of
characteristic zero. Suppose that every residue field of $\CC$ at a
closed point has cohomological dimension $\leq 1$. Then, $SK_1(C)$ is
divisible.
\end{proposition}

The proof of ($\cite{swan}$, Proposition 1.7) yields the following
result. 

\begin{proposition}\label{suslin}
Let $\CC = \Spec C$ be a curve as in (\ref{su1}). Then, the
natural homomorphism $K_1Sp(C) \ra SK_1(C)$ is an isomorphism. 
\end{proposition}


\section{Main Theorem}

Given an affine algebra $A$ over $\BR$
and a subset $I\subset A$, we denote by $Z(I)$, the closed subset of $X
= \Spec A$ defined by $I$ and by $Z_{\tr}(I)$, the set $Z(I)\cap X(\BR)$,
where $X(\BR)$ is the set of 
all real maximal ideals $\mm$ of $A$ (i.e. $A/\mm \iso \BR$).
We denote by $\Sing X$,  the set of all the prime
ideals $\p$ of $A$ such that $A_{\p}$ is not a regular ring.
\medskip

The following lemma is proved in ($\cite{O-P}$, Lemma 2).

\begin{lemma}\label{singular}
Let $A$ be a reduced affine algebra over $\BR$ of dimension $d$ and let 
$X=\Spec A$. Let
$u=(a_1,\ldots,a_n)$ be a unimodular row in $A^n$. Suppose $a_1 > 0$ 
on $X(\BR)$. Then, there exists $b_2,\ldots,b_n \in A$
such that $a_1+b_2a_2+\ldots+b_na_n > 0$ on $X(\BR)$
and $Z(a_1+b_2a_2+\ldots+b_na_n)$ is smooth on $X\backslash \Sing X$ of
dimension $\leq d-1$.
\end{lemma}

Now, we state the {\L}ojasiewicz's inequality ($\cite{BCR}$, 
Proposition 2.6.2).

\begin{lemma}\label{L-ineq}
Let $B$ be an affine algebra over $\BR$ and let $X=\Spec B$. Let
 $a,b\in B$ be such that $a/b$ is defined on a closed semi algebraic
 set $F\subset X(\BR)$. Then, there exists $g\in B$ such that $g >0$
 on $X(\BR)$ and $|a/b| < g$ on $F$.
\end{lemma}

The following lemma is an easy consequence of (\ref{serre}) and
(\ref{bass}). 

\begin{lemma}\label{lemma3}
Let $B$ be a ring of dimension $n$ and let $Q$ be a projective 
$B$-module of rank $n$. Let $J$ be an ideal of height $\geq
1$. Suppose $(a,q)\in \Um(B\op Q)$. Then, there exists $\Psi \in \Aut
(B\op Q)$ such that $(a,q) \Psi = (a_1,\wt q)$ with
$a_1=1$ (mod $J$) and $O(\wt q)=B$ (mod $J$).
\end{lemma}

\begin{proof}
Let ``bar'' denotes reduction modulo $J$. Since $\dim \ol B \leq n-1$ and
$\ol Q$ is a projective $\ol B$-module of rank $n$,
by Serre's theorem ($\ref{serre}$), $\ol Q$ has a unimodular element.
Let $\ol {q_1}\in \ol Q$ be a unimodular element, i.e.  
$O(\ol {q_1})=\ol B$. 

Since rank $\ol Q > \dim \ol B$, by Bass' theorem 
($\ref{bass}$), $E(\ol B\op \ol Q )$ acts
transitively on $\Um(\ol B \op \ol Q)$. Hence,
there exists $\ol \Psi \in E(\ol B \op \ol Q)$ such that 
$(\ol a ,\ol q)\ol \Psi = (1,\ol {q_1})$.

Applying ($\ref{one}$), $\ol \Psi$ can be lifted to an element
$\Psi \in \Aut (B \op Q)$.
 Let $(a,q)\Psi = (a_1,\wt q)$. Then, we have $a_1
= 1$ (mod $J$) and $O(\wt q)=B$ (mod $J$). This proves the result.
$\hfill \square$ 
\end{proof}

\begin{lemma}\label{lemma1}
Let $B$ be an affine algebra over $\BR$ and let 
$f\in B$ be an element not belonging to any real maximal ideal of $B$.
Let $K\subset B$ be an ideal and
$a\in B$ such that $f^r \in Ba+K$, for some integer $r$. Then, there
 exists $h\in 1+Bf$ such that $ah >0$ on $Z_{\tr}(K)$.
Moreover, if $I$ is any ideal of $B$ such that $f^l\in I+K$, for some
$l\in \BN$, then we can choose $h\in 1+If$.
\end{lemma}

\begin{proof}
Since $f^r\in Ba+K$, hence, $a$ has no zeros on $Z_{\tr}(K)$. Further,
it is given that $f^l \in
I+K$, hence, $f^{2l} = \lambda+\mu$, for some 
$\lambda\in I$ and $ \mu\in K$. Since $f^{2l} >0$ on $X(\BR)$,
where $X=\Spec B$, we get that the element $\lambda > 0$
on $Z_{\tr}(K)$. Applying lemma ($\ref{L-ineq}$)
for the element $1/af^2\gl$ (with $F=Z_{\tr}(K)$), we get an element 
$g\in B$ such that $g>0$ on $X(\BR)$ and $1/|a| f^2 \lambda  < g$ on
$Z_{\tr}(K)$. Thus, it follows that the element $(1+a f^2 \lambda g)a > 0$ 
on $Z_{\tr}(K)$. Let us write $h=1+af^2 \lambda g$. Then $h\in 1+If$.
Further, $ah > 0$ on $Z_{\tr}(K)$. This proves the lemma. 
$\hfill \square$
\end{proof}

\begin{lemma}\label{lemma2}
Let $B$ be an affine algebra over $\BR$ and let $X= \Spec B$. Let 
$f\in B$ be an element such that $f >0$ on $X(\BR)$. 
Let $K\subset B$ be an ideal and let
$a_1\in B$ be such that $a_1 >0$ on $Z_{\tr}(K)$. Then, there exists
$c\in K$ such that $a_1+c >0$ on $X(\BR)$. Moreover, if $J$ is any
ideal of $B$ such that $f^q - a_1 \in J$, for some $q\in \BN$, then we
can choose $c\in KJ$.
\end{lemma}

\begin{proof}
Let $W$ be the closed semi algebraic subset of $X(\BR)$
 defined by $a_1 \leq 0$. Let  $f^q-a_1= \nu \in J$. Since $f>0$ on
 $X(\BR)$, the element $\nu > 0$ 
on $W$. On the other hand, we have $Z_{\tr}(K)
\cap W = \varnothing$, since $a_1 \leq 0$ on $W$ and $a_1 >0$ on
 $Z_{\tr}(K)$. Hence,  if $K=(c_1,\ldots,c_n)$, then
$c_1^2+\ldots+c_n^2 > 0 $ on $W$. Therefore, applying ($\ref{L-ineq}$)
for the element  $a_1/\nu^2(c_1^2+ \ldots+c_n^2)$ (with $F=W$), we get 
an element $\wt c\in B$ such that $\wt c >0$ on $X(\BR)$ and 
$|a_1|/ \nu^2(c_1^2+\ldots+c_n^2) <\wt c$ on $W$. 
Let $c=\wt c \nu^2 (c_1^2+\ldots +c_n^2)$. Then $c\in KJ$. Further,
$a_1 + c > 0$ on $W$. 
We also have $a_1 + c>0$ on $X(\BR)\backslash W$,
 since $a_1>0$ on $X(\BR)/W$ and
$c \geq 0$ on $X(\BR)$. Therefore, we have $a_1+c > 0$ on the whole of
$X(\BR)$.
This proves the result.
$\hfill \square$
\end{proof}

\begin{lemma}\label{lemma4}
Let $B$ be an affine algebra over $\BR$ and let $I$ be an ideal of
$B$. Let $f\in B$ be an element such that $f>0$ on $X(\BR)$, where
$X=\Spec B$. Let $P$ be a projective $B_f$-module and let $(a,p)\in
\Um(B_{f} \op P)$ such that $a=1$ (mod $IB_f$) and $O(p)=B_f$ (mod
$IB_f$).  Then, there exists $h\in 1+Bf$ and $\Delta \in \Aut(B_{fh}
\op (P\ot B_{fh}))$ such that $(a,p)\Delta = (\wt a,p)$ with
$\wt a >0$ on $X(\BR)\cap \Spec B_{fh}$ and $\wt a=1$ (mod $IB_{fh}$).
\end{lemma}

\begin{proof}
Since $P$ is a projective $B_f$-module, we can find a $B$-module 
$M$ such that $P=M_f$.  Since
$(a,p)\in \Um(B_f\op M_f)$, 
after multiplying by a suitable power of $f$, we may assume that
$(a,p)\in B \op M$ such that $a = f^l$ (mod $IB$) and
$O(p) \supset  f^q B$ (mod $IB$), for some $l,q\in \BN$.

We have $(a,p)\in B\op M$ and $(a,p)_f \in \Um(B_f \op M_f)$.
Hence $f^r \in aB+O(p)$, for some $r\in \BN$. Write $K=O(p) B$.
We also have  $f^q\in K+I$. Hence, 
applying (\ref{lemma1}), there exists $h\in 1+fI$ such that $a_1
= h a > 0$ on $Z_{\tr}(K)$. 

Note that, we have $a_1 >0$ on $Z_{\tr}(K)$ and
$a_1 = f^l$ (mod $IB$) (since $h-1 \in I$).
Hence, applying (\ref{lemma2}), we get $c\in KI$ such that 
the element $a_2=a_1+c > 0$ on $X(\BR)$. Let $\gf\in P^*$ be such that
$\gf(p)=c$.
Note that, we still have
 $a_2 = f^l$ (mod $IB$). Let $\wt a=a_2/f^l \in B_f$. Then
$\wt a = 1$ (mod $IB_f$). 

From the above discussion, it is clear that
if $\Gamma_1=(h,Id),
\Gamma_2=(1/f^l,Id) \in \Aut (B_{fh}\op
P_{fh})$, then $(a,p)\Gamma_1=(a_1,p)$,
$(a_1,p)\Gamma_\gf = (a_2,p_1)$ and $(a_2,p)\Gamma_2 = (\wt
a,p)$. Further, $\wt a > 0$ on $X(\BR)\cap \Spec B_{fh}$ and
$\wt a = 1$ (mod $IB_f$). Take $\Delta=\Gamma_1
\Gamma_\gf\Gamma_2$. Then the result follows.
$\hfill \square$
\end{proof}
\medskip
 
The following result is an easy consequence of $(\ref{singular})$.

\begin{lemma}\label{lemma5}
Let $B$ be a reduced affine algebra of dimension $d$ over $\BR$ and let
$X=\Spec B$. Let $Q$ be a projective $B$-module. Let $J$ be the ideal
of $B$ defining the singular locus of $B$ and let $I\subset J$ be an ideal.
Let $(\wt a,q)\in \Um(B\op Q)$ such that $\wt a > 0$ on $X(\BR)$ and
$\wt a =1$ (mod $I$). Then, there exists $\Phi \in E(B\op Q)$ such that
$(\wt a,q)\Phi = (b,\wt q)$ with $b>0$ on $X(\BR)$, $Z(b)$ is smooth on
$X$ of dimension $\leq d-1$ and $\wt q\in IQ$.
\end{lemma}

\begin{proof}
Since $(\wt a,q)\in \Um(B\op Q)$, we have
$\wt a B + O(q)=B$. Further, $\wt a =1$ (mod $I$). Hence,
it is easy to see that if 
$I=(s_1,\ldots,s_l)$ and $O(q)=(c_1,\ldots,c_n)$, then
$(\wt a,s_1^2 c_1^2,\ldots,s_1^2 c_n^2,s_2^2
c_1^2,\ldots,s_l^2c_1^2,\ldots, s_l^2c_n^2)$
is a unimodular row in $B^{nl+1}$. 

Since $\wt a =1$ (mod $I$) and $I\subset J$, hence, $\wt a=1$ 
(mod $J$). Further, $\wt a >0$ on $X(\BR)$.
Applying ($\ref{singular}$), we get  $h_{ij}\in B$ such that 
$$b=\wt a + \sum_{i,j} h_{ij}s_i^2c_j^2 \,> 0$$ 
on  $X(\BR)$ and $Z(b)$ is smooth on $X$ of dimension $\leq d-1$.

Let $\gf\in Q^*$ be such that $\gf(q)=\sum_{i,j} h_{ij}s_i^2c_j^2$. Let
$\Delta_1=\Gamma_\gf \in E(B\op Q)$. Then
$(\wt a,q)\Delta_1 = (b,q)$. Note that $b=1$ (mod
$I$). Therefore, there
exists $\Delta_2 \in E(B\op Q)$ such that $(b,q)\Delta_2 =
(b,\wt q)$, where $\wt q \in IQ$. Write $\Phi=\Delta_1 \Delta_2$. Then
$(\wt a,q)\Phi = (b,\wt q)$ has the required properties. This proves
the lemma.
$\hfill \square$
\end{proof}
\medskip

The following result is a generalization of (\cite{O-P},  Proposition 1).

\begin{lemma}\label{positive}
Let $A$ be a reduced affine algebra of dimension $d$ over $\BR$ 
and let $X=\Spec A$. Let
$J$ be an ideal of $A$ of height $\geq 1$. Let $f\in A$ 
be an element not belonging to any real maximal ideal of $A$. 
Let $P$ be a projective 
$A_{f}$-module of rank $d-1$ and let $(a,p)\in \Um(A_{f} \op P)$.
Then, there exists $h\in 1+fA $ and $\Delta \in
\Aut(A_{fh} \op P_{fh})$ such that if 
$(a,p)\Delta = (b,\wt p)$, then 

$(1)$ $b > 0$ on $X(\BR)\cap \Spec A_{fh}$,

$(2)$ $Z(b)$ smooth on $\Spec A_{fh}$ of dimension $\leq
d-1$ and

$(3)$ $(b,\wt p)=(1,0)$ (mod $JA_{fh}$).
\end{lemma}

\begin{proof}
By replacing $f$ by $f^2$, we may assume that $f> 0$ on $X(\BR)$.
Let $J_1$ be the ideal of $A$ defining the singular locus of $\Spec
A$. Since $A$ is reduced and 
char $\BR=0$, $J_1$ is an ideal of height $\geq 1$. Let
 $I=JJ_1$. Then $\hh I \geq 1$.

Write $A_1=A_{f(1+fA)}$. Then $\dim A_1= d-1$. Recall, rank
$P=d-1$ and $(a,p)\in \Um(A_1\op (P\ot A_1))$.
Applying $(\ref{lemma3})$ with $B=A_1$, $Q=P\ot A_1$ and $J=IA_1$, 
there exists $\Psi \in \Aut (A_1 \op (P \ot A_1))$ such that
$(a,p)\Psi = (a_1,p_1)$, where $a_1
= 1$ (mod $IA_1$) and $O(p_1)=A_1$ (mod $IA_1$).
 
It is easy to see that there exists $h_1 \in 1+fA$ such that, if we
write $B=A_{h_1}$, then
$\Psi \in \Aut (B_f \op (P\ot B_f))$ and $(a,p)\Psi =
(a_1,p_1)$, where $a_1 = 1$ (mod $IB_f$) and
$O(p_1)=B_f$ (mod $IB_f$). Applying $(\ref{lemma4})$, there exists 
an element $h' \in 1+fB$ and $\Gamma \in \Aut(B_{fh'}\op (P\ot
B_{fh'}))$ such that $(a_1,p_1)\Gamma = (\wt a,p_1)$ with
$\wt a > 0$ on $X(\BR)\cap \Spec B_{fh'}$ and $\wt a = 1$ (mod $IB_{f}$)

Recall that $h_1\in 1+fA$, $B=A_{h_1}$ and $h'\in 1+fB$.
Let $s$ be an integer such that $h=h_1^s h'\in A$. Then $B_h=A_h$ and
$h\in 1+fA$. Therefore, there exists $\Gamma \in \Aut(A_{fh}\op (P\ot
A_{fh}))$ such that $(a_1,p_1)\Gamma = (\wt a,p_1)$ with
$\wt a > 0$ on $X(\BR)\cap \Spec A_{fh}$ and $\wt a = 1$ (mod $IA_{fh}$).

Write $C=A_{fh}$. Then, we have $(\wt a,p_1)\in \Um (C\op (P\ot C))$
such that $\wt a > 0$ on $X(\BR)\cap \Spec C$ and $\wt a =1$ (mod
$IC$), where $I\subset J_1$ (Recall that $J_1$ is an ideal of $A$
defining the singular locus of $\Spec A$). Applying $(\ref{lemma5})$,
we get $\Phi\in \Aut(C\op (P\ot C))$ such that 
$(\wt a,p_1) \Phi = 
(b,\wt p)$ with $b>0$ on $X(\BR)\cap \Spec C$, $Z(b)$ is smooth on
$\Spec C$ of dimension $\leq d-1$ and $\wt p \in I(P\ot C)$. 

Let $\Delta = \Psi\Gamma\Phi$.
Then $\Delta\in
\Aut(A_{fh}\op P_{fh})$ is such that $(a,p)\Delta=(b,\wt p)$ with $b >0$ on
$X(\BR)\cap \Spec A_{fh}$ and $Z(b)$ is smooth on $\Spec A_{fh}$
of dimension $\leq d-1$. Moreover, $(b,\wt p)=(1,0)$ (mod
$JA_{fh}$). This proves the lemma.  $\hfill
\square$
\end{proof}

\begin{proposition}\label{main}
Let $A$ be an affine algebra over $\BR$ of dimension $d$. Let 
$f\in A$ be an element not belonging to any real maximal ideal of $A$
and let $A' = A_{f(1+Af)}$. Then, every projective
$A'$-module $P$ of rank $d-1$ is cancellative.
\end{proposition}

\begin{proof}
By replacing $f$ by $f^2$, we may assume that $f> 0$ on $X(\BR)$,
where $X=\Spec A$.  It is enough to show that, if $(a,p)\in \Um(A'\op
P)$, then, there exists $\Lambda \in \Aut(A' \op P)$ such that
$(a,p)\Lambda = (1,0)$. Without loss of generality, we may assume that
$A$ is reduced.

We can choose $g \in 1+Af$ such that $P$ is a projective
$A_{fg}$-module of rank $d-1$. Write $\wt A=A_{fg}$. 
Let $t\in \wt A$ be a non-zero-divisor such that $P_t$ is a free $\wt
A_t$-module of rank
$d-1$. Let $F=\wt A^{d-1}$ be a free submodule of $P$ such that $F_t=P_t$
and let $s=t^l$ be such that $sP\subset F$. Let $(e_1,\ldots,e_{d-1})$
denote the standard basis of $\wt A^{d-1}$. 

Let $J$ be the ideal of $A$ defining the singular locus of $A$. Since
$A$ is reduced and char $\BR=0$, $J$ is an ideal of height $\geq 1$. 
Let $I=sJ$. Then $\hh I =1$. Applying 
(\ref{positive}) for $\wt B=A_g$, there exists $h'\in
1+f\wt B$ and $\Gamma \in \Aut (\wt B_{fh'}\op P_{fh'})$ such that 
if $(a,p)\Gamma = (a_1,p_1)$, then

(1) $a_1 > 0$ on $X(\BR)\cap \Spec \wt B_{fh'}$,

(2) $\Spec (\wt B_{fh'}/a_1 \wt B_{fh'})$ is smooth of dimension 
$\leq d-1$ and

(3) $(a_1,p_1)=(1,0)$ (mod $s \wt B_{fh'}$). \\

We can choose a suitable positive integer $r$ such that $h=g^r h' \in
A$ (in fact $h \in 1+Af$ and $\wt B_{fh}=A_{fh}$) and
$\Gamma \in \Aut (A_{fh}\op P_{fh})$ such that $(a,p)\Gamma = (a_1,p_1)$,  
satisfies the above properties $(1)-(3)$, i.e.

(1) $a_1 > 0$ on $ X(\BR)\cap \Spec A_{fh}$,

(2) $\Spec (A_{fh}/a_1 A_{fh})$ is smooth of dimension $\leq d-1$ and

(3) $(a_1,p_1)=(1,0)$ (mod $sA_{fh}$).  \\

Since $p_1\in sP_{fh}$ and we have $sP \subset F$, hence, we can write
$p_1 = b_1e_1+\ldots+ b_{d-1}e_{d-1}$, for some $b_i\in A_{fh}$. 
Let $B=\BR(f)\otimes_{\tr[f]} A_{fh}$. Then $B$ is an affine algebra over
$\BR(f)$ of dimension $d-1$. We write $\wt P = P\ot B$.

Let ``bar'' denotes reduction modulo the ideal $a_1B$. Since $a_1B+sB=B$ and
$F_s=P_s$, it follows that
the inclusion $F\subset P$ gives rise to the equality $\ol
F = \ol P$. In particular, $\ol P$ is free of rank $d-1$ with a basis
$(\ol e_1,\ldots,\ol e_{d-1})$ and $\ol p_1$ is a unimodular element
of $\ol P$.

If $d=3$, then  $C=B/a_1B$ is smooth of dimension $1$. Note that, 
every maximal
ideal $\mm$ of $C$ is the image in $\Spec C$ of a prime ideal
$\p$ of $A_{fh}$ of height $2$ containing $a_1$. Since $a_1$ does not
belongs to any real maximal ideal of 
$A_{fh}$, by (\cite{serre}), the
residue field $\BR(\p)=k(\mm)$ has cohomological dimension $\leq 1$.
By (\ref{suslin}), $SK_1(C)$ is divisible and the natural
map $K_1Sp(C) \ra SK_1(C)$ is an isomorphism.
Hence, there exists $\Gamma' \in \SL_2(C) \cap ESp_4(C)$ and
$c_1,c_2\in B$ such that, if $ q=c_1^2e_1 + c_2e_2\in F$, then
$\Gamma'(\ol p_1) = \ol q$. By ($\ref{two}$),
$\Gamma'$ has a lift $\Gamma_1\in SL(\wt P)$. Recall that $\wt P= P\ot
B$. \\

If $d \geq 4$, then, since $\ol P$ is a free of rank $d-1$, using
($\ref{su}$) and (\ref{suslin}), one can deduce from the proof 
of ($\cite{Su2}$, Theorem
2.4) that there exists $\wt \Gamma \in E(\ol P)$ and $c_i \in
B$, $1\leq i\leq d-1$ such that, if $q=c_1^{d-1}e_1+c_2e_2+\ldots+
c_{d-1}e_{d-1}\in F$, then $\wt\Gamma(\ol p_1) = \ol q$.
By ($\ref{one}$),
 $\wt \Gamma$ can be lifted to an element
$\Gamma_1 \in SL(\wt P)$. (In particular, the above argument shows that 
every stably free $B/a_1B$-module of rank $\geq d-2$ is cancellative).

Therefore, in either case, there exists $q_1\in \wt P$ and $\Gamma_1
\in \SL(\wt P)$ such that
$$\Gamma_1(p_1) =q-a_1q_1, \;{\rm where}\; q=c_1^{d-1}e_1 
+c_2e_2 + \ldots+c_{d-1}e_{d-1}.$$

Now, the rest of the argument is similar to (\cite{Bh}, Theorem 4.1).
We give the proof for the sake of completeness.

Now, $\Gamma_1$ induces an automorphism $\Psi_1 = (Id_B,\Gamma_1)$ of $B\op
\wt P$. Let $\wt \Psi = \Psi_1\,\Delta_{q_1}$. Recall that $\Delta_{q_1}\in
E(B\op \wt P)$ is defined as $\Delta_{q_1}(b,q')=(b,q'+bq_1)$.
Therefore, we have $(a_1,p_1)\wt \Psi = (a_1,q)$. Let us write
$\Lambda_1=\Gamma \wt \Psi \in \Aut(B\op \wt P)$ . Then $(a,p) \Lambda_1
=(a_1,q)$.

Recall that $a_1=1$ (mod $sB$). Hence, there exists 
$x\in B$ such that $sx+a_1=1$. Let
$\mu_i =sxc_i$. Then $\mu_i-c_i\in a_1B$. Let 
$$q_2
= \mu_1^{d-1}e_1+ \sum_2^{d-1}\mu_ie_i \in sF~~~ \rm{and}~~~ 
q_3 = \sum_1^{d-1}\mu_i e_i.$$ 
Then $q_2-q = a_1 p_2$, for some $p_2\in F$.
Hence, we have $(a_1,q)\Delta_{p_2} = (a_1,q_2)$. Let $\Lambda_2 =
\Lambda_1 \Delta_{p_2}$. Then $(a,p)\Lambda_2 = (a_1,q_2)$.

Since $1-a_1\in sB$ and $\mu_i \in sB$ for $1\leq i\leq d-1$. 
Hence, the row $[a_1,\mu_1,\ldots,\mu_{d-1}]
\in \Um_d(B,Bs)$. Therefore, by (\ref{kallen}),
$$[a_1^{d-1},\mu_1,\ldots,\mu_{d-1}] =
[a_1,\mu_1^{d-1},\ldots,\mu_{d-1}]\;\; (\rm{mod}\;\, E_d(B,Bs)).$$ 

By ($\ref{lift}$), 
there exists $\Phi \in E(B\op \wt P)$ such
that $(a_1,q_2)\Phi = (a_1^{d-1},q_3)$. Write $\wt \Phi = \Lambda_2
\Phi$. Then, we have
$(a,p) \wt \Phi = (a_1^{d-1},q_3)$.

Since $\wt P/a_1\wt P$ is free of rank $d-1$ and every stably free
$B/a_1B$-module of rank $\geq d-2$ is cancellative, $\ol q_3\in
\wt P/a_1 \wt P$ can be extended to a basis of $\wt P/a_1 \wt P
$. Therefore, by ($\ref{su}$), 
there exists $\Phi_1 \in \Aut (B\op
\wt P)$ such that $(a_1^{d-1},q_3) \Phi_1=(1,0)$. 

Let $\Lambda =\wt \Phi\Phi_1.$  Then $\Lambda\in \Aut(B\op
\wt P)$ and $(a,p)\Lambda=(1,0)$. Note that $A'=A_{f(1+Af)} = B\ot_{R(f)}
A'$. Therefore, we get the result.
$\hfill \square$        
\end{proof}
\medskip

As a consequence of the above proposition (\ref{main}), we prove the
following result. If  $P=A^{d-1}$ in the following theorem (\ref{cor}), then
we get ($\cite{O-P}$, Theorem).

\begin{theorem}\label{cor}
Let $A$ be an affine algebra over
$\BR$ of dimension $d$ and let $f\in A$ be an element not belonging to
any real maximal ideal of $A$. Let $P$ be a projective $A$-module
of rank $\geq d-1$. Let $(a,p)\in A_{f}\op P_f$ be a
unimodular element. Then, the projective $A_{f}$-module $Q=A_{f}\op
P_f/(a,p)A_{f}$ is extended from $A$.
\end{theorem}

\begin{proof}
Let $A'=A_{f(1+Af)}$.  By (\ref{main}), $P\ot A' \iso Q\ot A'$. Hence, there
exists $g\in 1+Af$ and an isomorphism $\Psi : P\ot A_{fg} \iso Q\ot
A_{fg}$.  The
module $Q$ over $A_f$ and $P$ over $A_g$ together with an isomorphism
$\Psi$ yield a projective module over $A$ whose extension to $A_f$
is isomorphic to $Q$. This proves the result.  $\hfill \square$.
\end{proof}

\remark The above theorem ($\ref{cor}$) is valid for an affine algebra
$A$ over any real closed field $k$. For simplicity, we have taken
$k=\BR$.     \\

\noindent{\bf Acknowledgments.}
I sincearly thank Prof. S. M. Bhatwadekar for suggesting the problem and
for useful discussion.

{}


\begin{thebibliography}{}

\bibitem{B}{} H. Bass, {\it K-Theory and stable algebra},
I.H.E.S. {\bf 22} (1964), 5-60.

\bibitem{BCR}{} J. Bochnak, M. Coste, M.-F. Roy, {\it Real algebraic
geometry}, Ergebnisse der Mathematik, Berlin Heidelberg New York,
Springer, 1998.
 
\bibitem{Bh}{} S. M. Bhatwadekar, {\it A cancellation theorem for
projective modules over affine algebras over $C_1$-fields}, 
J.P.A.A. {\bf 183} (2003), 17-26.

\bibitem{Bh-Roy}{} S. M. Bhatwadekar and A. Roy, {\it Some theorems
about projective modules over polynomial rings}, J. Algebra {\bf 86}
(1984), 150-158.

\bibitem{K}{} W. van der Kallen, {\it A group structure on certain
orbit sets of unimodular rows}, J. Algebra {\bf 82} (1983), 363-397.
 
\bibitem{O-P}{} M. Ojanguren and R. Parimala, {\it Projective modules
over real affine algebras}, Math. Ann. {\bf 287} (1990), 181-184.

\bibitem{S1}{} J. P. Serre, {\it Sur les modules projectifs}, Sem.
Dubreil-Pisot {\bf 14} (1960-61), 1-16.

\bibitem{serre}{} J. P. Serre, {\it Sur la dimension cohomologique des
groupes profinis}, Topology {\bf 3} (1968), 264-277.

\bibitem{su1}{} A. A. Suslin, {\it A cancellation theorem for
projective modules over affine algebras}, Sov. Math. Dokl. {\bf 18}
(1977), 1281-1284.

\bibitem{Su2}{} A. A. Suslin, {\it Cancellation over affine
varieties}, J. of Soviet Math. {\bf 27} (1984), 2974-2980.

\bibitem{swan}{} R. G. Swan, {\it A cancellation theorem for
projective modules in the metastable range}, Invent. Math. {\bf 27}
(1974), 23-43.

\bibitem{V}{} L. N. Vaserstein, {\it Operation on orbit of unimodular
vectors}, J. Algebra {\bf 100} (1986), 456-461. 
\end{thebibliography}
\end{document}